\documentclass[]{cuparticle}
\makeatletter\if@twocolumn\PassOptionsToPackage{switch}{lineno}\else\fi\makeatother

\usepackage{tabulary}
\usepackage{amsmath}
\usepackage{amssymb}
\usepackage{amsfonts}
\usepackage{fontenc}
\usepackage[utf8]{inputenc}
\usepackage{graphicx}
\usepackage{epstopdf}
\usepackage{url}
\usepackage{multirow}
\usepackage{morefloats}
\usepackage{floatflt}
\usepackage{cancel}
\usepackage{tfrupee}
\makeatletter

\AtBeginDocument{\@ifpackageloaded{textcomp}{}{\usepackage{textcomp}}}
\makeatother
\usepackage{colortbl}
\usepackage{xcolor}
\usepackage{pifont}
\usepackage[nointegrals]{wasysym}
\makeatletter

\def\mcWidth#1{\csname TY@F#1\endcsname+\tabcolsep}

\def\cAlignHack{\rightskip\@flushglue\leftskip\@flushglue\parindent\z@\parfillskip\z@skip}
\def\rAlignHack{\rightskip\z@skip\leftskip\@flushglue \parindent\z@\parfillskip\z@skip}

\usepackage{ifxetex}
\ifxetex\else\if@twocolumn\usepackage{dblfloatfix}\fi\fi

\AtBeginDocument{
\expandafter\ifx\csname eqalign\endcsname\relax
\def\eqalign#1{\null\vcenter{\def\\{\cr}\openup\jot\m@th
  \ialign{\strut$\displaystyle{##}$\hfil&$\displaystyle{{}##}$\hfil
      \crcr#1\crcr}}\,}
\fi
}

\AtBeginDocument{%
  \@ifpackageloaded{endfloat}%
   {\renewcommand\efloat@iwrite[1]{\immediate\expandafter\protected@write\csname efloat@post#1\endcsname{}}}{\newif\ifefloat@tables}%
}%

\def\BreakURLText#1{\@tfor\brk@tempa:=#1\do{\brk@tempa\hskip0pt}}
\let\lt=<
\let\gt=>
\def\processVert{\ifmmode|\else\textbar\fi}

\@ifundefined{subparagraph}{
\def\subparagraph{\@startsection{paragraph}{5}{2\parindent}{0ex plus 0.1ex minus 0.1ex}%
{0ex}{\normalfont\small\itshape}}%
}{}

\newcommand\role[1]{\unskip}
\newcommand\aucollab[1]{\unskip}

\@ifundefined{tsGraphicsScaleX}{\gdef\tsGraphicsScaleX{1}}{}
\@ifundefined{tsGraphicsScaleY}{\gdef\tsGraphicsScaleY{.9}}{}
\def\checkGraphicsWidth{\ifdim\Gin@nat@width>\linewidth
	\tsGraphicsScaleX\linewidth\else\Gin@nat@width\fi}

\def\checkGraphicsHeight{\ifdim\Gin@nat@height>.9\textheight
	\tsGraphicsScaleY\textheight\else\Gin@nat@height\fi}

\def\fixFloatSize#1{}
\let\ts@includegraphics\includegraphics

\def\inlinegraphic[#1]#2{{\edef\@tempa{#1}\edef\baseline@shift{\ifx\@tempa\@empty0\else#1\fi}\edef\tempZ{\the\numexpr(\numexpr(\baseline@shift*\f@size/100))}\protect\raisebox{\tempZ pt}{\ts@includegraphics{#2}}}}

\AtBeginDocument{\def\includegraphics{\@ifnextchar[{\ts@includegraphics}{\ts@includegraphics[width=\checkGraphicsWidth,height=\checkGraphicsHeight,keepaspectratio]}}}

\DeclareMathAlphabet{\mathpzc}{OT1}{pzc}{m}{it}

\def\URL#1#2{\@ifundefined{href}{#2}{\href{#1}{#2}}}

\def\UrlOrds{\do\*\do\-\do\~\do\'\do\"\do\-}%
\g@addto@macro{\UrlBreaks}{\UrlOrds}

\edef\fntEncoding{\f@encoding}

\makeatother

\makeatletter
\gdef\@recto{}
\gdef\@verso{}
\@typesetRunningHeads{}
\def\historyattrib#1{}
\makeatother

\volume{}
\doi{}

\begin{document}

\authorheadline{Yi Wang}
\runningtitle{G-networks and the optimization of supply chains}

\begin{frontmatter}

\title{G-networks and the optimization of supply chains}
\author[1]{Yi Wang}

\address[1]{Intelligent Systems and Networks Group,\\
Department of Electrical and Electronic Engineering,\\
Imperial College, London SW7 2AZ, UK\\
\texttt{yi.wang18@imperial.ac.uk}}


\begin{abstract}
Supply chains are fundamental to the economy of the world and many supply chains focus on perishable items, such as food, or even clothing that is subject to a limited shelf life due to fashion and seasonable effects. G-networks have not been previously applied to this important area. Thus in this paper, we apply G-networks to supply chain systems and investigate an optimal order allocation problem for a N-node supply chain with perishable products that share the same order source of fresh products. The objective is to find an optimal order allocation strategy to minimize the purchase price per object from the perspective of the customers. An analytical solution based on G-networks with batch removal, together with optimization methods are shown to produce the desired results. The results are illustrated by a numerical example with realistic parameters.
\end{abstract}

\end{frontmatter}
\providecommand{\keywords}[1]
{
  \small	
  \textbf{\textit{Keywords---}} #1
}
\\

\section{Introduction}
A supply chain is a system of organizations, people, activities, information, and resources involved in moving a product or service from suppliers to customers \cite{frazzon2015synchronizing}, and in \cite{braziotis2013supply} the distinction between supply chains and supply networks is clarified. Traditional supply chains focus on inventory control strategies, system profit and customer demand \cite{li2008simulating}. Simulation based optimization is often suggested to determine the safety stock level to meet customer satisfaction \cite{jung2004simulation}. In \cite{yin2015optimal}, one manufacturer and multiple suppliers are studied to seek an optimal discount policy for coordinating the relationship between manufacturers and suppliers. The closed-loop supply chains include traditional forward supply chains and reverse supply chains and many researchers have addressed the importance and potential of this area \cite{guide2003challenge}. The management and optimization of two competitive closed-loop supply chains are discussed in \cite{gu2012management} and in \cite{yang2018product} an exploratory case study is proposed to highlight the need of practical tools to guide practitioners in identifying value opportunities in circular cycles.

In recent years, the strategic importance of perishable goods has led to the development of perishable product supply chains \cite{duan2013new}. Thus, in \cite{kopach2008tutorial}, an age-based policy using queuing theory is suggested to construct a red blood cell inventory system with urgent demand rates and non-urgent demand rates. In \cite{sarker2000supply}, a mathematical model is developed to determine an ordering policy for deteriorating products under inflation, permissible delay of payment and allowable shortage. In \cite{haijema2009blood}, a combined stochastic dynamic programming-simulation approach is suggested to research the production and inventory management of blood platelets under irregular production breaks. In \cite{duan2013new}, a policy based on old inventory ratio is developed to minimize the outdate rate under a maximal allowable shortage level. In \cite{rong2011optimization}, a food supply chain model is used to plan food production and distribution via mixed-integer linear programming and in \cite{besik2017quality} a path-based framework for food supply chains based on kinetics and game theory is considered. Time critical supply chains in the Australian dairy industry are studied in \cite{ivanov2016dynamic} and in \cite{chen2013performance} a supply chain simulation study for agricultural distribution systems is examined.

The importance of information and communication systems can help support greater adaptation and resilience in current supply chains \cite{bremang2006information}. In \cite{ivanov2012inter} an inter-disciplinary modelling framework for multi-structural and collaborative cyber-physical networks is discussed and in \cite{xu2009modelling} the product information tracking and feedback via wireless technologies is modelled. Moreover, many researchers have presented the potential of information technologies such as RFID system and CPS \cite{nemeth2006adopting}, \cite{frazzon2015synchronizing}, \cite{puschmann2005successful}. Note that one of the key features of a supply chain is that multiple products can share multiple facilities with capacity constraints and demands from multiple customers \cite{jung2004simulation}, which is overlooked by most research. The development of information and communication systems provides convenience for the demand distribution and the facility sharing of supply chains.

Furthermore, high information flow and product flow among different supply chains lead to the development of supply chain risk management (SCRM). In \cite{rajesh2015modeling}, supply chain risk management is viewed as a blend of grey relational analysis and decision-making and in \cite{hulsmann2008strategic} the potential of autonomous cooperation and control (ACC) is discussed to strike a balance between flexibility and stability for global supply chains (GSC). In \cite{saenz2014creating}, an example of Cisco Systems is given to present the importance of risk management for companies to create more resilient supply chains and in \cite{rajesh2015selection} the major risks in electronic supply chains and risk mitigation strategies are discussed via digraph-matrix approaches. However, most research focuses on conceptual aspects and lacks effective mathematical models and algorithms. Besides, few researchers consider related optimization problems about the message security of supply chains, which has become more important with the progress of information technologies \cite{Cyber1}.

In this paper, the class of queueing networks knows as G-networks with batch removal are introduced into the area of supply chains and a multiple-node system is described to show how to model supply chains. An optimization problem of multiple-node supply chains is detailed with perishable products and the sharing of a same order source. The target of this problem is to minimize the buying price per unit object and an analytical solution is obtained via Lagrange multipliers techniques.

The remainder of this paper is organized as follows. In Section 2, we summarize the basic results and applications of G-networks and present a multiple-node supply chain system. In sections 3, we detail an optimization problem via G-networks with batch-removal. We solve the optimization problem by Lagrange multipliers, obtain an analytical solution and illustrate a numerical example in Section 4. Conclusions are given in Section 5.

\section{The G-Network Model}
A G-network is an open network of $N$ queues with independent and exponentially distributed
(i.i.d.) service times at each queue with rate $\mu(i)\geq 0$ for the $i-th$ queue, in which customers are either ``positive'' or ordinary customers which request and receive service, or they are ``signals'' \cite{GN2}. Signals may be negative customers, triggers, resets or adders, as discussed in a series of papers \cite{GN4,GN5,gelenbe2002g,GN6}. In this paper we only consider those signals which are triggers or which are negative customers with batch removal. They were first introduced as a model of spiking neuronal networks \cite{Stable}.

A positive customer can add the queue length by one and a negative customer can reduce the queue length by one through destroying a positive customer or forcing it to leave the network, if the queue length is positive. After finishing service, a positive customer from queue $i$ may head for queue $j$ as a positive customer with probability $p^{+}(i,j)$, or as a signal with probability $p^{-}(i,j)$, or it will depart from the network with probability $d(i)$. The movement of customers in the network is described by a Markov chain. In the G-Network, external arrivals  occur according to Poisson processes. Again, external arrivals
are either positive customers arriving at a queue $i$ according to a Poisson process at rate $\Lambda(i)$, or they are signals that
arrive according to a Poisson process at rate $\lambda(i)$.

The simplest result concerning a G-Network that only has positive and negative customers can be stated as follows \cite{GN1}.
Let $K(t)=(K_1(t),~...~,K_N(t)$ denote the number of positive customers at the $N$ queues of the G-Network, and let
$k=(k_1,~...~,k_N)$ be a particular value that {\it i} may take, with the $k_i\geq 0$. Also, define the stationary probability distribution
$p(k) =\lim_{t\rightarrow\infty}Prob [K(t) = k]$ if it exists.

\medskip

\noindent{\bf Result 1} Consider the following system of non-linear equations:
\begin{eqnarray}
&&\Lambda^{+}(i) = \Lambda(i) + \sum_{j=1}^{N}{q_{j} \mu_{j} p^{+}(j,i)},\\
&&\lambda^{-}(i) = \lambda(i) + \sum_{j=1}^{N}{q_{j} \mu_{j} p^{-}(j, i)},\\
&&where~q_i=\frac{\Lambda^+(i)}{\mu_i+\Lambda^-(i)},~i\in\{1,~...~,N\}.
\end{eqnarray}
If the above system of equations has a solution such that all the $0\leq q_i<1$ then:
\begin{equation}
p(k) = \prod_{i=1}^{n}[1-q_{i}]q_{i}^{k_{i}}~.
\end{equation}
Further details and the proof can be found in \cite{GN1,GN3,Gelenbe-Pujolle}.

\medskip

Further features have been introduced into G-networks such as:
\begin{itemize}
\item Triggers which can move positive customers from queue $i$ to some other queue $j$ or force them to leave the network \cite{GN4},
\item Batch removal which represents situations when one ``negative customer''  can force a batch of customers of random size to leave the network from a given queue \cite{GN5},
\item Adders which can probabilistically change the queue length at the service center that it visits acting as a load supplier or regulator \cite{GN6}, and
\item Resets which can reset the queue length to a random queue length whose distribution is identical to the stationary distribution at that queue if the queue is empty \cite{gelenbe2002g}.
\item  In \cite{gelenbe1998g}, the concept of multiple classes of negative and positive customers has been introduced into G-Networks, together with three types of service centers and service disciplines (first-in-first-out, processor sharing and last-in-first-out).
\end{itemize}

G-networks have had many applications. A few prominent ones include the random neural network (RNN) which is a spiking ``integrate and fire'' model for combinatorial optimisation, video compression web search, deep learning and image processing \cite{Combinatorial,VideoCompression,Cloud,web,gelenbe2016deep}, the energy packet network (EPN) which can be used to provide energy on demand to the cloud computing servers \cite{EPN,NOLTA,Fourneau}, and the cognitive packet network (CPN) which makes use of adaptive techniques and machine learning to seek out routes based on defined QoS criteria in cyber-physical systems and computer networks \cite{gelenbe2002cognitive,SAN,Cloud}, and can also be used to minimize energy consumption in networks \cite{gelenbe2011energy}.

G-Networks have also been used to model adversarial agents in applications such as Cybsersecurity \cite{gelenbe2007dealing} and Gene Regulatory Networks \cite{GelenbePhysRev2010}.

\subsection{Multiple-node supply chains}

In this part, we consider a general case of supply chain systems containing $N$ nodes with one class of order messages and one class of objects, which may be perishable products or normal products. A node $i$ can receive objects from some other node $j$ with a probability $p^{+}(j,i)$, where $i,j=1,...,N$, or from outside of the supply chain with a Poisson rate $\Lambda(i)$. Besides, the node can receive an order message from some other node $j$ with a probability $p^{-}(j,i)$ or from outside of the supply chain at a Poisson rate $\lambda(i)$. Every order message triggers a batch of objects from a node $i$ to some other node $j$ if the arrival queue is not empty. The order message disappears and has no effect if it is sent into an empty queue.

A node $i$ can send an object or an order message to other node $j$ at a service rate $\mu(i)$ with a probability $p^{+}(i,j)$ or a probability $p^{-}(i,j)$, where $\sum_{j=1}^{N}{p^{+}(i,j)}+{p^{-}(i,j)}=1$. The service rate $\mu(i)$ obeys independent and exponential distribution $(i.i.d.)$. The perishable rate of node $i$ is $\beta_i$, which means the object is wasted because of deterioration, expiry or other reasons. Note that perishable products have a larger perishable rate than normal products.

As for an order message, the batch size $B(i)$ is determined by the order demand of customers. The distribution of batch size $B(i)$ is given by $P[B(i) = s] = \pi_{i}$, where $s \geq 1$. If the queue length of objects in this node is $k(i)$ and $k(i) \geq B(i)$, then its length is reduced by $B(i)$; if $k(i) \leq B(i)$, then the queue length becomes zero. We use two nodes to show the object traffic and the message traffic in this multiple-node supply chain, which is presented in Figure 1.
\begin{figure}
\centerline{\includegraphics[height=7cm,width= 10cm]{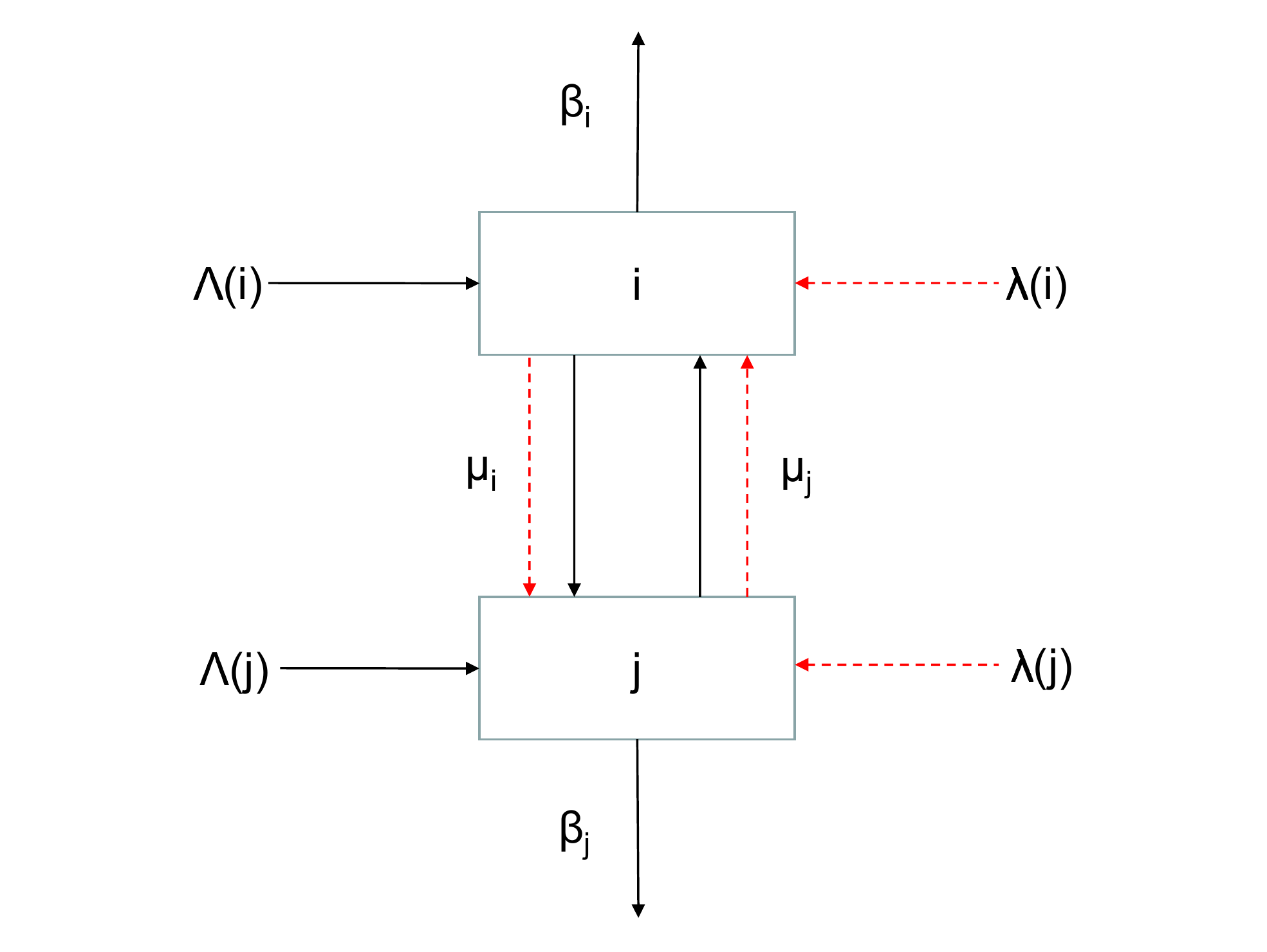}}
\caption{The simple presentation of a multiple-node supply chain with two nodes $i$ and $j$. Red lines represent message traffic and black lines represent object traffic. Note that in this general case node $i$ sends objects to node $j$ at a rate $\mu_i p^+(i,j)$ and node $i$ sends messages to node $j$ at a rate $\mu_i p^-(i,j)$. Objects arrive to node $i$ at a rate of $\Lambda_i$ and the arrival rate of messages in node $i$ is $\lambda_i$. The perishable rate of node $i$ is $\beta_i$.}
\end{figure}

In order to analyze this system, we use a basic result concerning G-network models with batch removal \cite{GN5} which generalizes Result 1:

\medskip

\noindent{\bf Result 2} Consider the following system of non-linear equations:
\begin{eqnarray}
&&\Lambda^{+}(i) = \Lambda(i) + \sum_{j=1}^{N}{q_{j} \mu_{j} p^{+}(j,i)},\\
&&\lambda^{-}(i) = \lambda(i) + \sum_{j=1}^{N}{q_{j} \mu_{j} p^{-}(j, i)},\\
&&where~q_i=\frac{\Lambda^+(i)}{\mu_i + \beta_i +\Lambda^-(i) \frac{1-F(q)}{1-q}},~i\in\{1,~...~,N\},\\
&&and~F(q)=\sum_{r=1}^\infty \pi(r) q^r~. \label{F(q)}
\end{eqnarray}
If the above system of equations has a solution such that all all the $0\leq q_i<1$ then:
\begin{equation}
p(k) = \prod_{i=1}^{N}[1-q_{i}]q_{i}^{k_{i}}~.
\end{equation}
Note that the effect of batch removal appears in (\ref{F(q)}) where $F(q)$ is the generating function for the size of the batch of objects which are removed after each sale.

\section{The optimization problem of a multiple-node supply chain}
Some of the above considerations can be included in a mathematical model which offers the possibility to optimize the system by choosing appropriate parameters. In this section, we consider a $N$ node supply chain controlled by a single retailer and each node represents a warehouse. Note that there are no movements of objects or messages among these warehouses and they share a same order source.

The objects are received from a manufacturing facility or from some suppliers, and the arrival process of objects to a warehouse $i$ is assumed to be Poisson at rate $\lambda(i)$. Orders to purchase objects arrive to warehouses $i$ in the form of a Poisson process at rate $\gamma P_i$, where $\sum_{i=1}^N P_i = 1$, and each order is a request to buy a batch of objects of size $R$, which is assumed to be a random variable with probability distribution $P[R=r]=\pi(r)$, $\sum_{r=1}^\infty \pi(r)=1$ and finite average: $E[R]<\infty$. We assume that successive orders are i.i.d..

Obviously, it may not be possible to satisfy all of the requests of a given order simply because the warehouse may not contain enough objects when a purchase order arrives. Thus we call $S_i$ the number of objects which are actually bought when a purchase order arrives.

In the case studied in this section, we assume objects are perishable products, which are wasted with an exponentially distributed ``service'' time at rate $\mu_i$. The batches (triggered by $\gamma$) refer to batch orders of goods which are sold for a price $C$, where the price $C$ is a fixed price $c_0$ plus a variable price related to the number of objects which are actually bought. The structure of this optimization problem is presented in Figure 2 and the steady probability function related to $q_i$ is given by:
\begin{eqnarray}
&&q_i = \frac{\lambda_i}{\mu_i + \gamma P_i \frac{1-F(q_i)}{1-q_i}},\\
&&where~F(q_i)=\sum_{r=1}^\infty \pi(r) \cdot q_i^r~. \label{F(q)}
\end{eqnarray}
A customer expects to buy more objects using a lower price. In other words, the customer hopes to minimize the buying cost per unit object. Clearly, the buying cost per unit object will be:
\begin{eqnarray}
&&C = c_0 + a\sum_{i=1}^{N} \gamma P_i E[S_i],\\
&&where~a<0.
\end{eqnarray}
where $E[S_i]$ is the expected number of objects taken from warehouse $i$ after each order, $c_0$ is the initial buying price per unit object and $a$ is the cost coefficient related to the total number of objects taken by customers. The parameter $a$ should be negative to represent that a customer can get a lower price if he buys more objects.
\begin{figure}
\centerline{\includegraphics[height=7cm,width= 10cm]{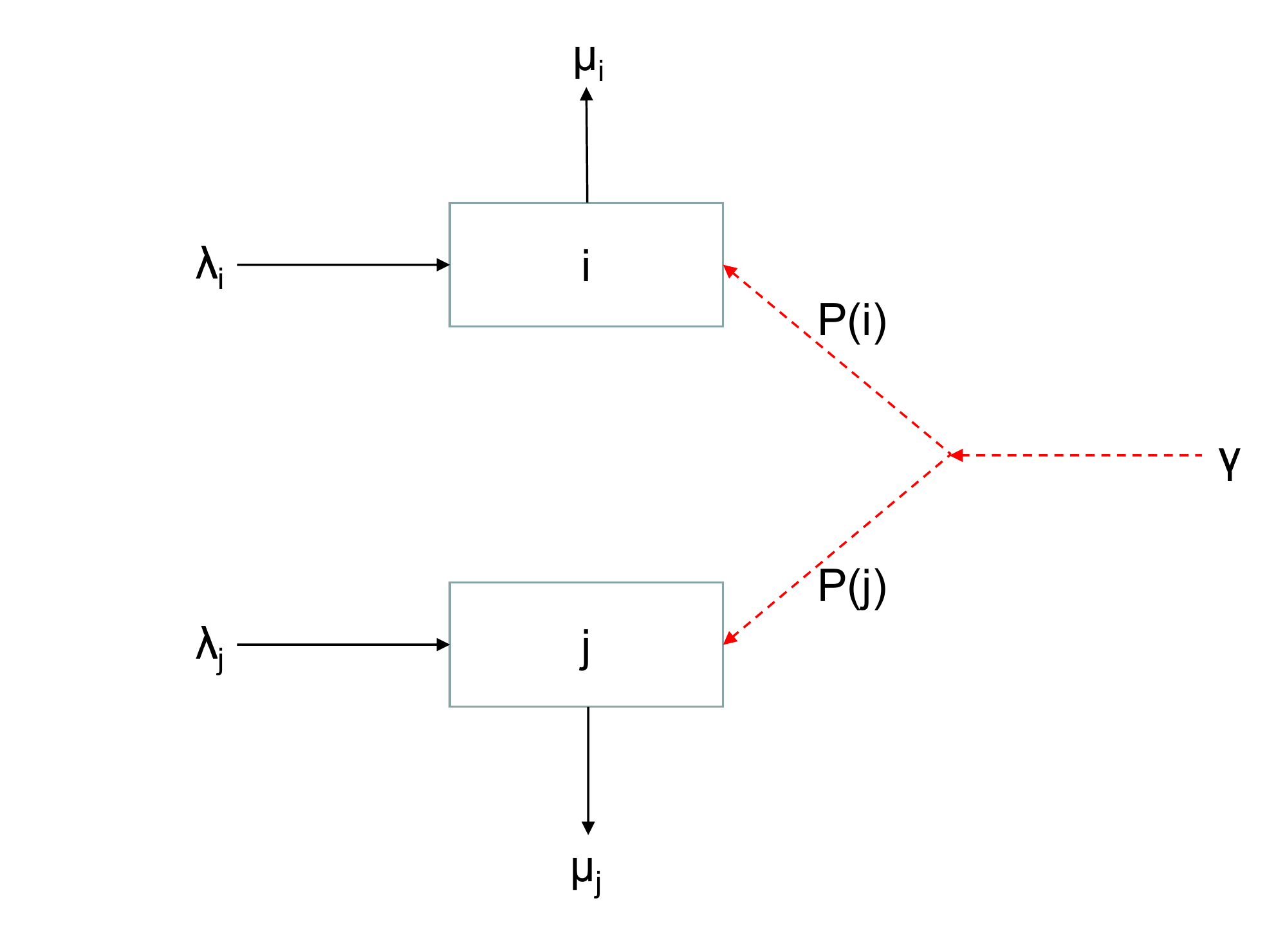}}
\caption{A schematic representation of the optimization problem with two warehouses $i$ and $j$. Red lines represent order traffic and black lines represent object traffic. The order arrival rate at node $i$ is $\gamma P_i$ and objects arrive to node $i$ at a rate of $\lambda_i$. Objects in node $i$ have a perishable rate of $\mu_i$. Note that we assume there are no movements of objects or messages between these two warehouses.}
\end{figure}

Note that if a customer orders $R$ objects, she/he can only buy $K$ objects if $K\leq R$. On the other hand, the customer can order $R$ objects and buy them all when $K>R$. Generalising the result of single-node supply chain in \cite{gelenbesupply}, we have:
\begin{eqnarray}
E[S_i~|~R=r]&=&\sum_{k=0}^{r}kq_i^k(1-q_i)+r\sum_{k=r+1}^\infty q_i^k(1-q_i),\nonumber\\
&=&\sum_{k=0}^{\infty}kq_i^k(1-q_i)-\sum_{k=r+1}^\infty kq_i^k(1-q_i)+rq_i^{r+1},\\
&=&q_i\frac{1-q_i^r}{1-q_i}.
\end{eqnarray}
Hence by using the definition of $F(q_i)=\sum_{k=1}^\infty \pi(r)q_i^r$, the expected number of objects bought by an order can be computed as:
\begin{eqnarray}
&&E[S_i]=q_i\frac{1-F(q_i)}{1-q_i}.
\end{eqnarray}

\section{Analysis of the optimization problem}
The above optimization problem considers $N$ warehouses sharing a single order source and the target is to minimize the buying cost per unit object with the constraint $\sum_{i=1}^NP_i=1$. In other words, the target is to maximize the number of objects taken by customers, because the unit price decreases with the increase of the number of bought objects. To analyse this problem, we assume that the number of required objects per unit order obeys a geometric distribution $\pi(r)=u^r\frac{1-u}{u}$ where $0<u<1$, and remember that $r=1,~2,~...$. Using formula (10), (\ref{F(q)}) and (12), we have:
\begin{eqnarray}
&&F(q_i)=\frac{q_i(1-u)}{1-q_iu},\\
&&q_i=\frac{\lambda_i}{\mu_i + \gamma P_i \frac{1}{1-q_iu}},\\
&&C= c_0 + a \sum_{i=1}^{N}\gamma P_i E[S_i],\\
&&~~= c_0 + a \sum_{i=1}^{N}(\lambda_i-q_i \mu_i),\\
&&where~a<0.
\end{eqnarray}
However, using formula 18, we obtain:
\begin{eqnarray}
&&q_i=\frac{\lambda_i}{\mu_i+\frac{\gamma P_i}{1-q_iu}},\\
&&~~=\frac{\gamma P_i+\mu_i+\lambda_i u}{2u\mu_i} \pm \sqrt{(\frac{\gamma P_i+\mu_i+\lambda_i u}{2u\mu_i})^{^2}-\frac{\lambda_i}{u\mu_i}}~.
\end{eqnarray}
Because of the existence condition of $q_i$ ($0\leq q_i<1$) for $i=1,...,N$, we have:
\begin{eqnarray}
&&q_i=\frac{\gamma P_i+\mu_i+\lambda_i u}{2u\mu_i} - \sqrt{(\frac{\gamma P_i+\mu_i+\lambda_i u}{2u\mu_i})^{^2}-\frac{\lambda_i}{u\mu_i}}~,\\
&&~~=A_i-\sqrt{A_i^2-\frac{\lambda_i}{\mu_iu}},\\
&&where~A_i=\frac{\gamma P_i+\mu_i+\lambda_i u}{2u\mu_i}\label{A_i1}.
\end{eqnarray}
Following the above calculations, we use Lagrange multiplier techniques to solve this optimization problem. Let a real number $\beta$ be the Lagrange multiplier and the modified cost function based on lagrangian is written as:
\begin{eqnarray}
&&\psi = C + \beta(\sum_{i=1}^NP_i-1),\\
&&~~= c_0 + a \sum_{i=1}^{N}(\lambda_i-q_i \mu_i)+\beta(\sum_{i=1}^NP_i-1).
\end{eqnarray}
Suppose $P^*=(P_1^*,...P_N^*)$ is a local solution of this problem, the necessary Kuhn-Tucker conditions are:
\begin{eqnarray}
&&\nabla_p\psi(P^*,\beta^*)=0,\\
&&\sum_{i=1}^N{P_i^*-1}=0.
\end{eqnarray}
The partial derivatives of $\psi$ and $q_i$ with respect to any $P_i$, where $i=1,...,N$, is given by:
\begin{eqnarray}
&&\frac{\partial \psi}{\partial P_i} = -a \mu_i q_i^{P_i}+\beta,\label{C}\\
&&\frac{\partial q_i}{\partial P_i}=\frac{\gamma}{2\mu_iu}(1-\frac{A_i}{\sqrt{A_i^2-\frac{\lambda_i}{\mu_iu}}})\label{q_i}.
\end{eqnarray}
Substitute formula (\ref{C}) and (\ref{q_i}) into the necessary Kuhn-Tucker conditions, we get the following formulas:
\begin{eqnarray}
&&-a \mu_i \frac{\gamma}{2\mu_iu}(1-\frac{A_i}{\sqrt{A_i^2-\frac{\lambda_i}{\mu_iu}}})+\beta=0,\\
&&(B^2-1)A_i^2=\frac{\lambda_i}{\mu_iu}B^2,\\
&&A_i=\sqrt{\frac{\lambda_i}{\mu_iu}}\sqrt{1+\frac{1}{B^2-1}}~,\label{A_i2}\\
&&where~B=1-\frac{2u\beta}{a \gamma}.
\end{eqnarray}
Substitute formula (\ref{A_i1}) into formula (\ref{A_i2}), the solution $P_i^*$ is:
\begin{eqnarray}
&&P_i^*=\frac{2 \mu_i u \sqrt{\frac{\lambda_i}{\mu_iu}}\sqrt{1+\frac{1}{B^2-1}}-\mu_i-\lambda_iu}{\gamma}~.
\end{eqnarray}
Moreover, the second necessary condition $\sum_{i=1}^N{P_i^*-1}=0$ is presented by:
\begin{eqnarray}
\sum_{i=1}^N \frac{2 \sqrt{\lambda_i \mu_i u}\sqrt{1+\frac{1}{B^2-1}}-\mu_i-\lambda_iu}{\gamma}=1,\\
\sqrt{1+\frac{1}{B^2-1}}=\frac{\gamma+\sum_{i=1}^N(\mu_i+\lambda_iu)}{\sum_{i=1}^N 2 \sqrt{\lambda_i \mu_i u}}\label{B}.
\end{eqnarray}
Substitute formula (\ref{B}) into formula (38), we obtain:

\medskip

\noindent{\bf Result 3}
The optimal solution to this problem is given by:
\begin{eqnarray}
&&P_i^*=\frac{\sqrt{\lambda_i \mu_i}}{\gamma}~[\frac{\gamma+\sum_{i=1}^N(\mu_i+\lambda_iu)}{\sum_{i=1}^N \sqrt{\lambda_i \mu_i}}]-\frac{\mu_i+\lambda_iu}{\gamma}.
\end{eqnarray}
However, the sufficient condition that exists an optimum solution $p^{*}$ also needs to be examined. To guarantee the existence of the local minimum, the Hessian matrix $\nabla_{PP}\psi$ must be positive definite. Note that $\nabla_{PP}\psi$ is a diagonal matrix with diagonal entries:
\begin{eqnarray}
&&\frac{\partial^2\psi(P^*,\beta^*)}{\partial P_i^2}=\frac{\partial^2 C}{\partial P_i^2}=\frac{-a \gamma^2 \lambda_i}{4 \mu_i^2 u^3 (A_i^2-\frac{\lambda_i}{\mu_i u})^{\frac{3}{2}}}.
\end{eqnarray}
Because $a<0$, the sufficient condition holds if the following inequality is satisfied for $i=1,...,N$:
\begin{eqnarray}
&&A_i^2>\frac{\lambda_i}{\mu_i u}.
\end{eqnarray}
According to the existence of $q_i$, the inequality is strictly satisfied and the existence of the strict constrained local minimum is proved. Note that $q_i < 1$, using formula (25) and (26), we obtain:
\begin{eqnarray}
&&P_i>\frac{(\lambda_i-\mu_i)(1-u)}{\gamma}.
\end{eqnarray}
This condition is physically meaningful since it implies that the demand for objects is sufficiently large to consume objects in these warehouses to keep this supply chain stable.

\subsection{A Numerical example}
In order to illustrate the analytically obtained optimal solution of this optimization problem, a numerical example considering three warehouses is presented and the related parameters are shown in Table 1. In this case, we assume the products are perishable objects, which have different perishable rates in different warehouses. The node 1 represents the biggest warehouse with a big arrival rate of objects and and a low perishable rate, while the node 3 represents the smallest one with a small arrival rate of objects and a high perishable rate. Note that these three warehouses share a same arrival rate of orders and this numerical example uses a sufficiently large $\gamma$ to guarantee the existence of $q_i$, where $0\leq q_i<1$ and $i=1,2,3$.
\begin{table}
\centerline{Tab. 1. Parameters of the supply chain system with three nodes}
\centerline{\includegraphics[height=4cm,width= 6cm]{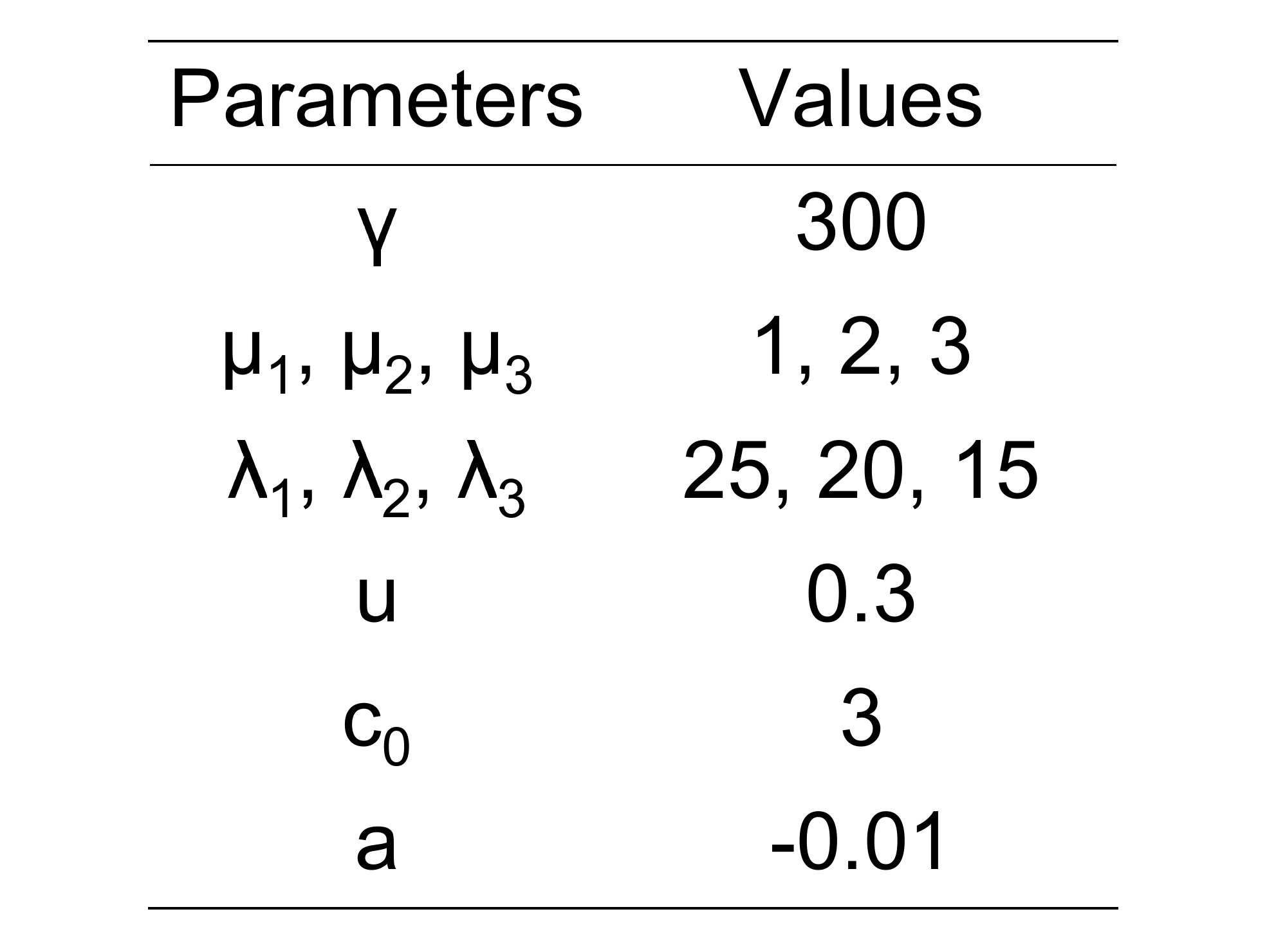}}
\end{table}

Using the constraint $P_1+P_2+P_3=1$ and Result 3, the optimal solution is $P^*=(P_1^*,P_2^*,P_3^*)=(0.2711,0.3521,0.3768)$ and the final price per unit object is $2.4100$, compared with the initial buying price per unit object $c_0=3$. Using formula (43), we have numerical conditions to guarantee the existence of $q_i$ ($0.056<P_1<1$, $0.042<P_2<1$, $0.028<P_3<1$).The result is presented in Figure 3 where the x-axis and y-axis are $P_1$ and $P_2$ separately and the z-axis is the buying cost per unit object $C$. To make it more clear, the neighbourhood of the optimum point $P^*$ is also presented in Figure 3, which illustrates our analytical solution gets the optimal point to minimize the buying cost per unit object. Besides, the influences of $\lambda_i$ and $\mu_i$ on $P_i$, where $i=1,2,3$, are showed in Figure $4$. We can find $P_i$ increases when $\lambda_i$ and $\mu_i$ become bigger and the increasing trend becomes slower with the change of $\lambda_i$ and $\mu_i$.
\begin{figure}
\centering
\begin{tabular}{cc}
\begin{minipage}[t]{0.5\textwidth}
\centering
\includegraphics[height=4.5cm,width=6cm]{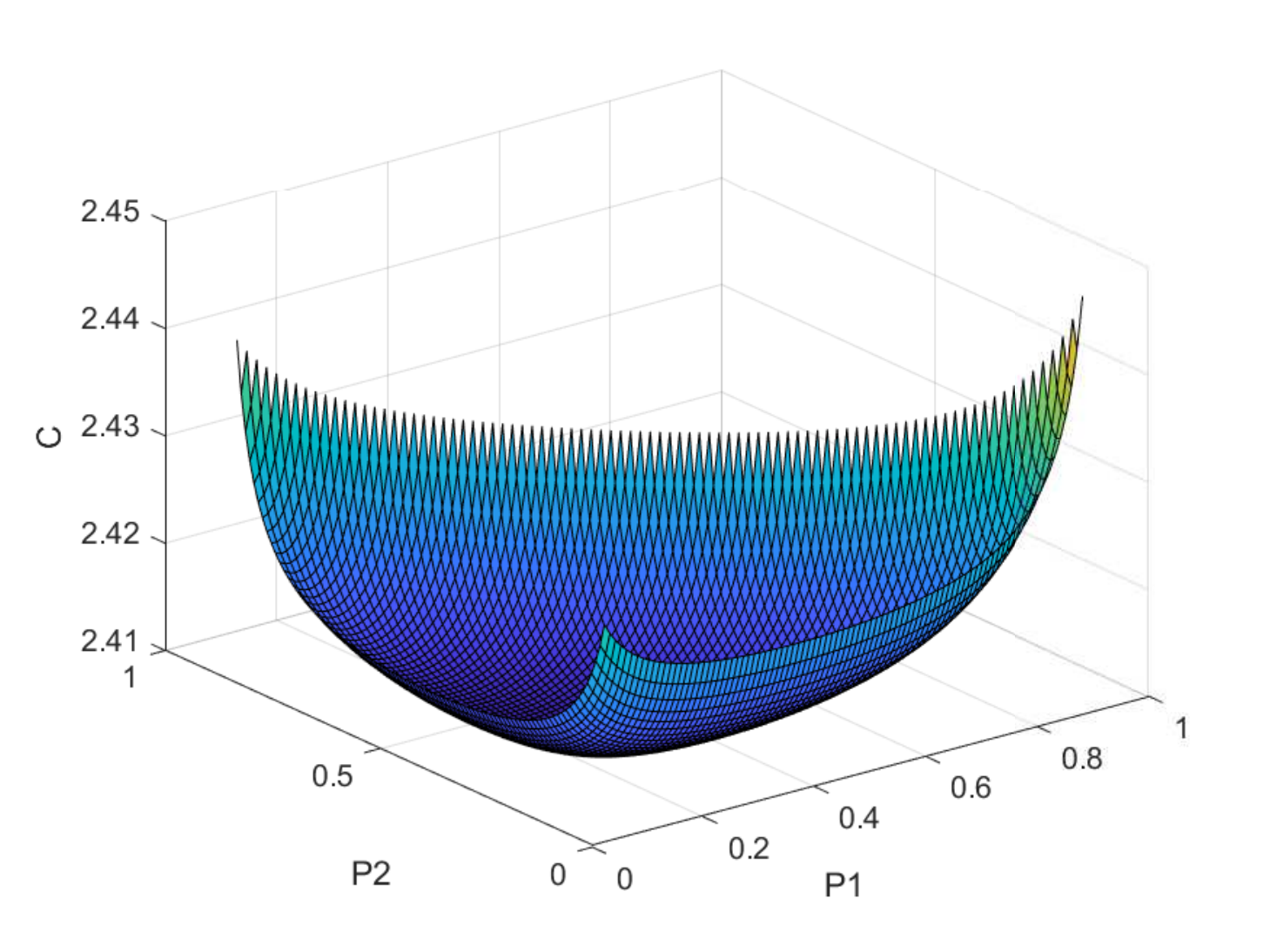}
\end{minipage}
\begin{minipage}[t]{0.5\textwidth}
\centering
\includegraphics[height=4.5cm,width=6cm]{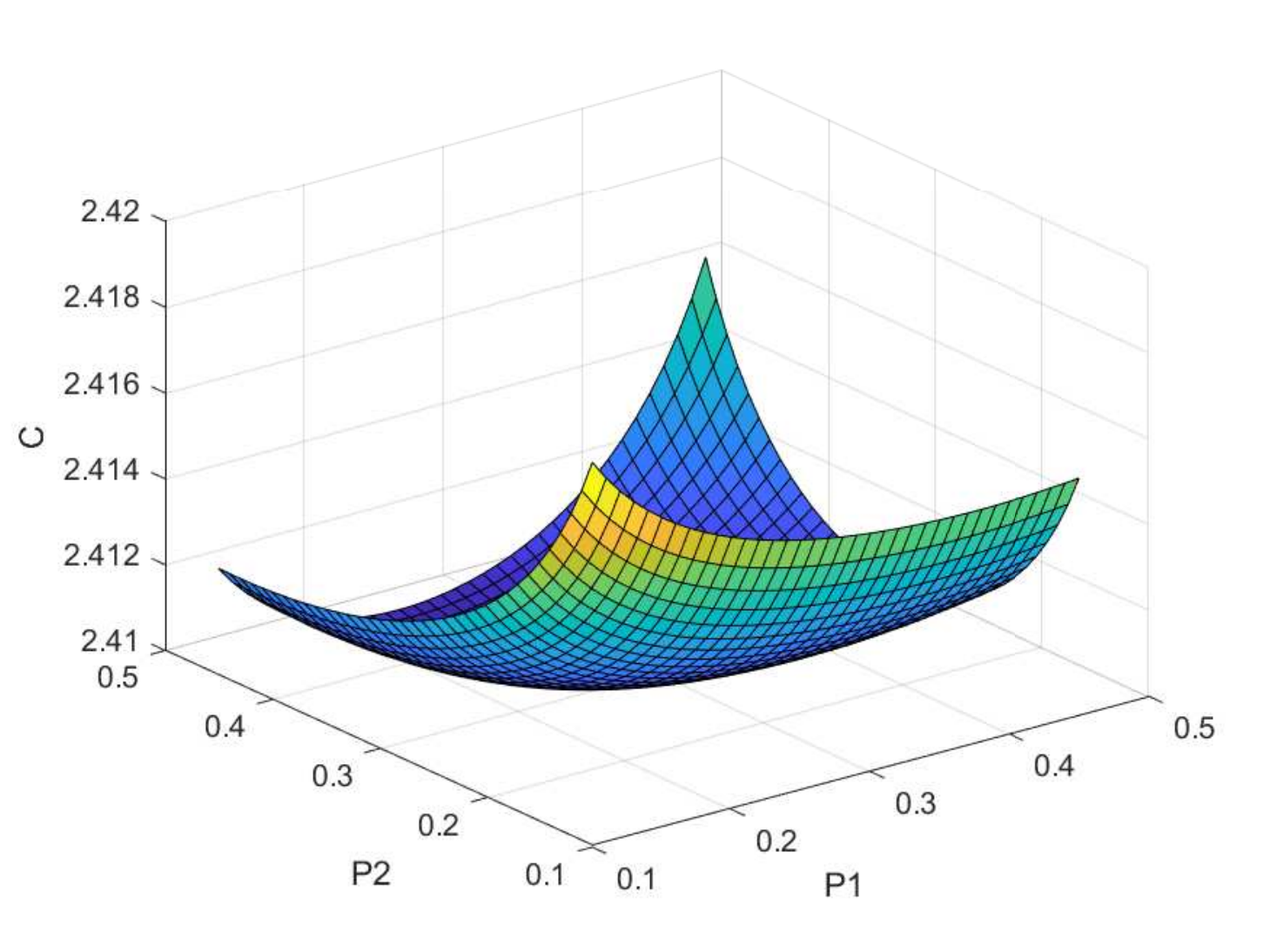}
\end{minipage}
\end{tabular}
\caption{The buying cost per unit object C for all $(P_1,P_2)$ pairs. In the left figure, the $(P_1,P_2)$ pairs can get values from 0 to 1 under the constraint $P_1+P_2+P_3=1$. In the right figure, the $(P_1,P_2)$ pairs get values in the range of 0.1 to 0.45 to illustrate the existence of the optimum point $P^*$ more clearly.}
\end{figure}

\begin{figure}
\centering
\begin{tabular}{cc}
\begin{minipage}[t]{0.5\textwidth}
\centering
\includegraphics[height=4.5cm,width=6cm]{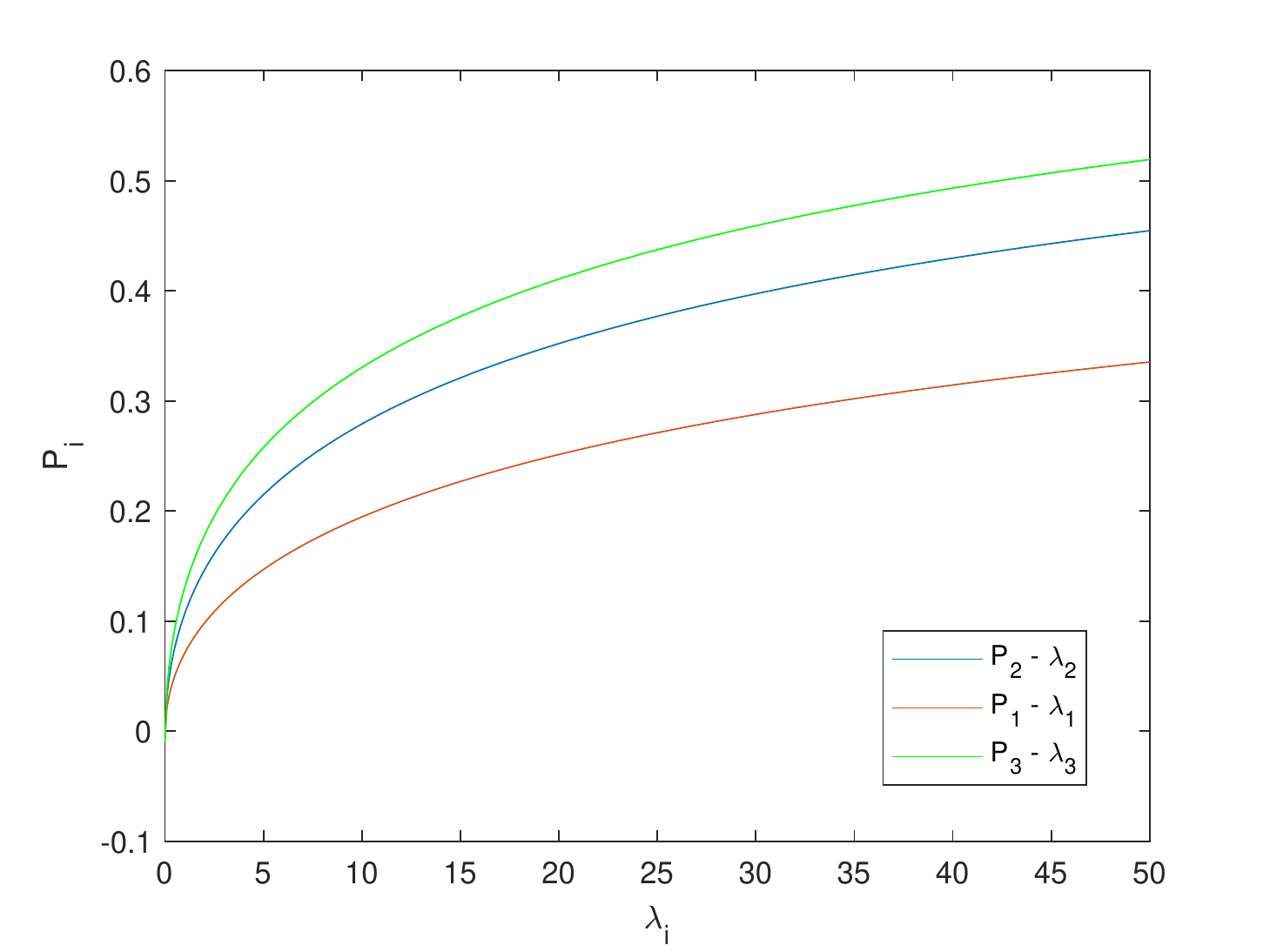}
\end{minipage}
\begin{minipage}[t]{0.5\textwidth}
\centering
\includegraphics[height=4.5cm,width=6cm]{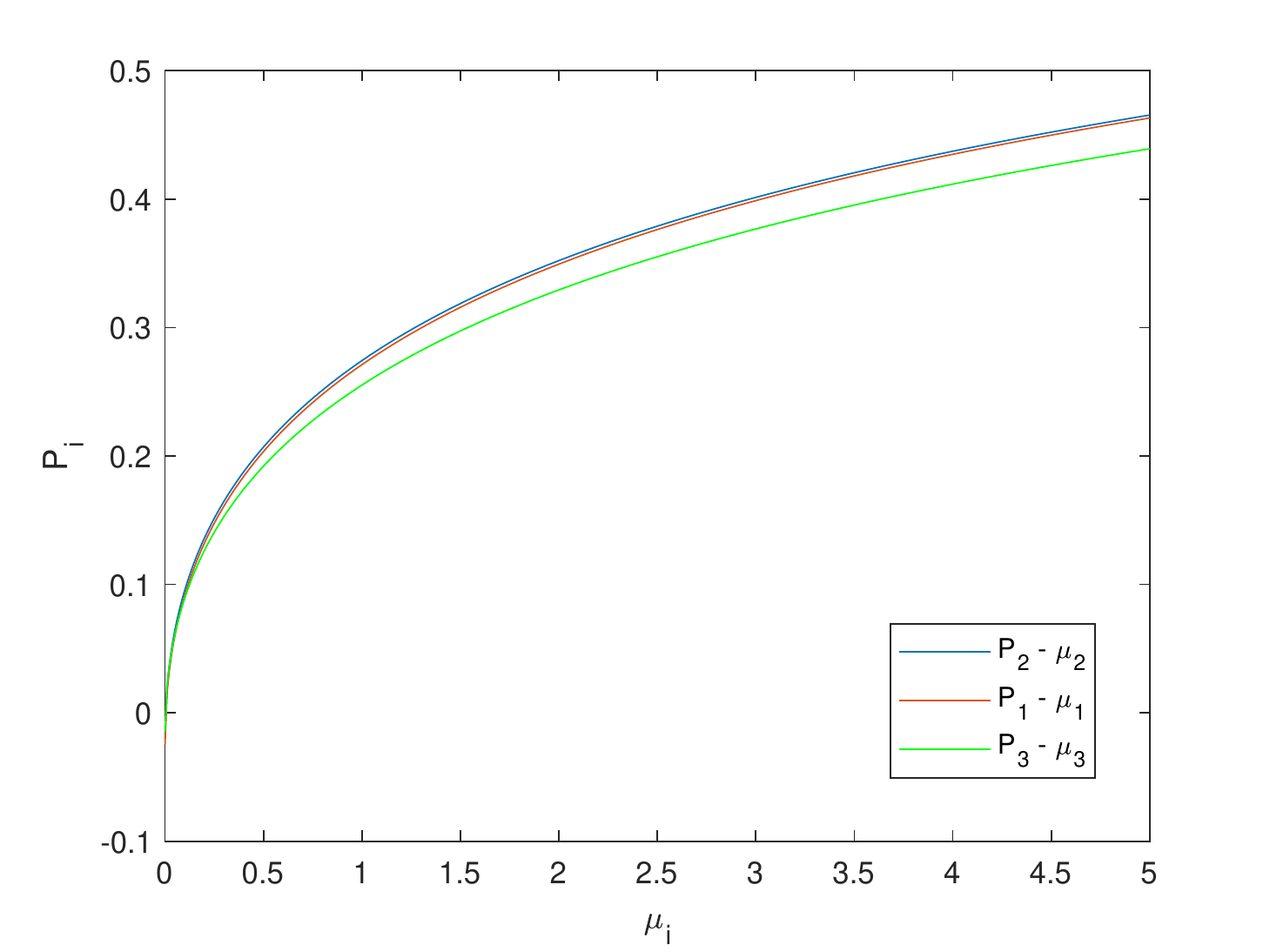}
\end{minipage}
\end{tabular}
\caption{The influence of $\lambda_i$ and $\mu_i$ on $P_i$, where $i=1,2,3$. The left figure shows $P_i$ increases when $\lambda_i$ becomes bigger, which shows a bigger arrival rate of objects may attract more orders. The right figure shows $P_i$ also increases with the increase of $\mu_i$, which means a customer tends to send more orders to a warehouse with a relatively high perishable rate in order to buy more objects.}
\end{figure}

Furthermore, the influences of $P_i$ on $q_i$ and $E[S_i]$ is presented in Figure 5 where $i=1,2,3$. The values of $q_i$ and $E[S_i]$ decrease with the increase of $P_i$. Besides, there is $q_1=0.2785$, $q_2=0.1762$ and $q_3=0.1246$, when $P^*=(P_1^*,P_2^*,P_3^*)=(0.2711,0.3521,0.3768)$. Under this situation, the customer buys 24.7215 objects per unit time from warehouse 1, 19.6477 objects per unit time from warehouse 2 and 14.6263 objects per unit time from warehouse 3 averagely.

\begin{figure}[htbp]
\centering
\begin{tabular}{cc}
\begin{minipage}[t]{0.5\textwidth}
\centering
\includegraphics[height=4.5cm,width=6cm]{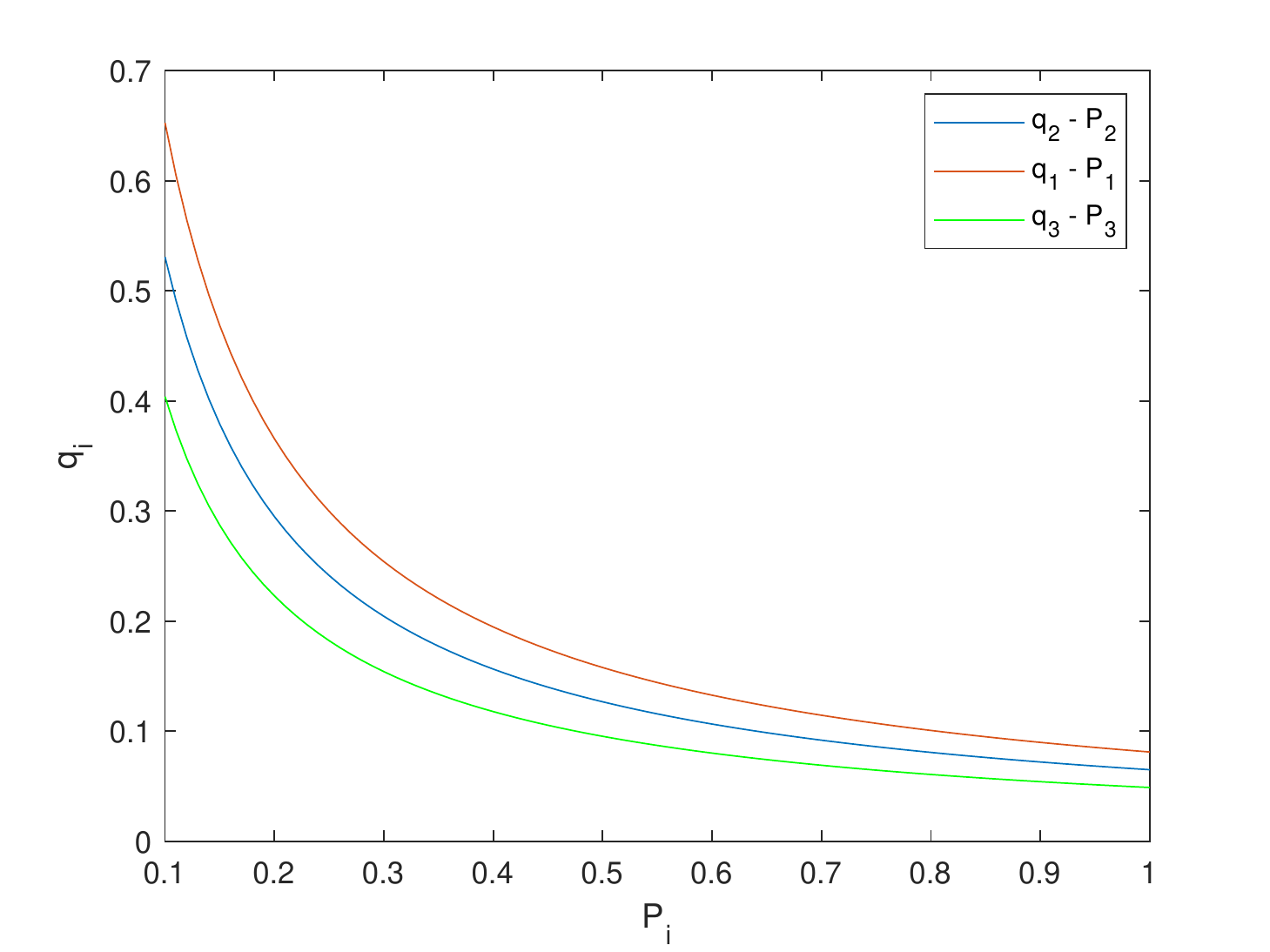}
\end{minipage}
\begin{minipage}[t]{0.5\textwidth}
\centering
\includegraphics[height=4.5cm,width=6cm]{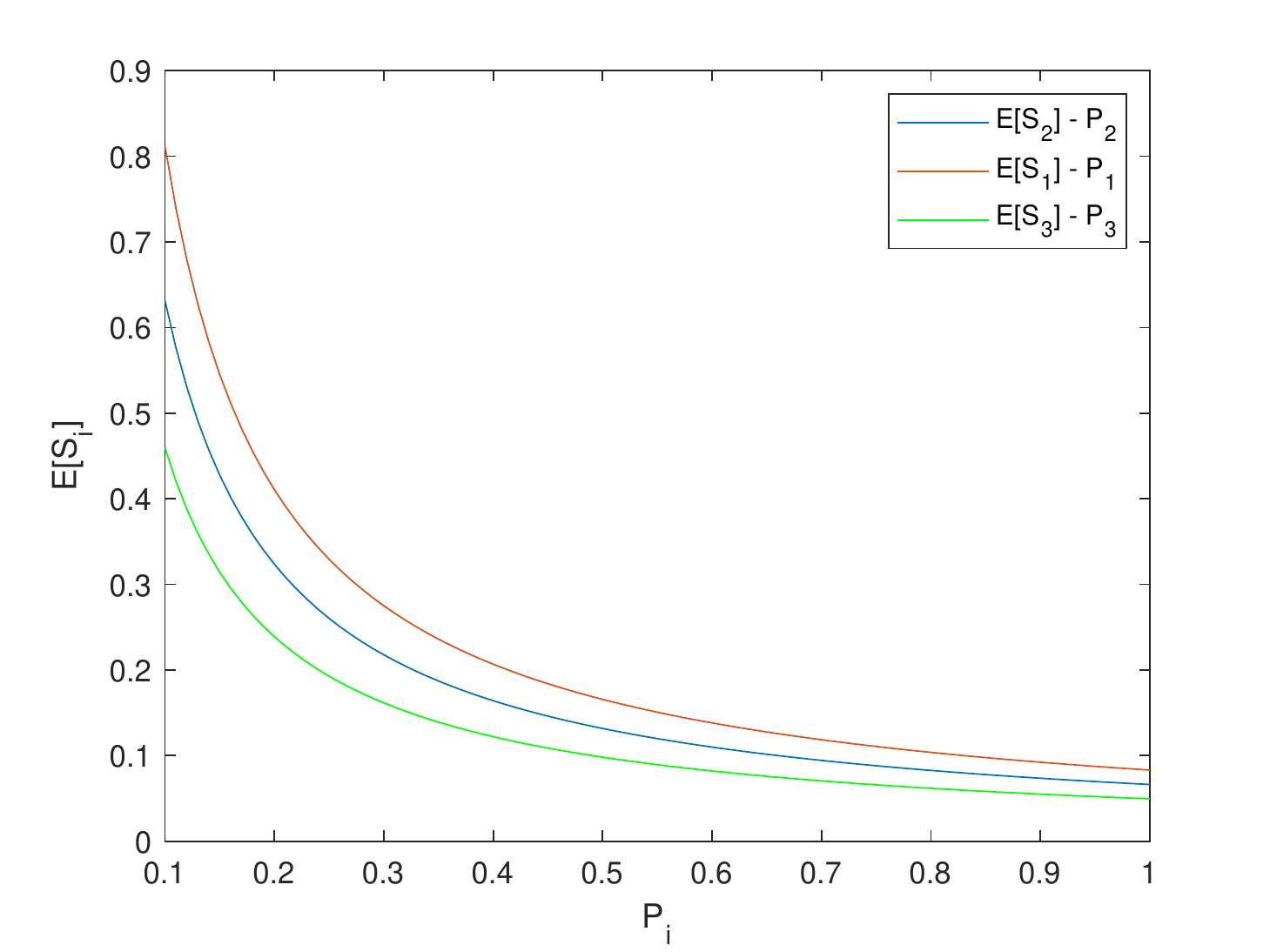}
\end{minipage}
\end{tabular}
\caption{The influence of $P_i$ on $q_i$ and $E[S_i]$, where $i=1,2,3$. The left figure shows $q_i$ decreases when $P_i$ increases and the right figure shows $E[S_i]$ will decrease when $P_i$ increases. These two figures show a bigger $P_i$ takes more objects per unit time from node $i$.}
\end{figure}

\section{Conclusions}

Supply chains are special instances of Cyber-Physical Systems \cite{Wu} which are devoted to the forwarding and commercialisation of various forms of goods. In this paper, we have briefly reviewed the literature on supply chains and we apply G-networks with batch removal to model and optimize such systems. We proposed a mathematical model based on G-networks to represent certain aspects of such systems, focusing on an optimization problem of multiple-node supply chains considering perishable products and the sharing of a same order source and show how to model the supply chain system with G-networks to minimize the buying price per unit object from the view of customers. Note that in this case, we assume objects can not be moved from one warehouse to another warehouse and we only consider one type of objects. We use Lagrange multipliers techniques to solve this problem and get an analytical solution of $P^*$. We also produce a numerical example with three retail nodes to illustrate our results.

Future work will investigate optimization problems considering multiple classes of nodes and multiple classes of objects and study how to guarantee a profit per unit time. The security problem caused by message traffic of supply chains should also be considered in the future \cite{gelenbe2007dealing}. Moreover, future work should also model the connections and data interdependencies among different types of nodes \cite{Sevcik} to optimize supply chain systems from an overall systemic perspective.

\section*{Acknowledgement}

The author thanks Prof. E. Gelenbe for his help in formulating and solving the problems discussed in this paper.

\bibliographystyle{unsrt}

\end{document}